\documentclass[journal]{IEEEtran}

\usepackage{cite}

\ifCLASSINFOpdf
\else
\fi
\usepackage{array}

\usepackage{amssymb,amsmath}
\usepackage{graphicx}
\usepackage{amsthm}
\usepackage{epstopdf}
\allowdisplaybreaks[4]
\usepackage{flushend}

\emergencystretch=\maxdimen
\newtheorem{theorem}{Theorem}

\newtheorem{lemma}{Lemma}

\newtheorem{assumption}{Assumption}
\newtheorem{remark}{Remark}

\newcommand{\col}{$\upshape{col}$}
\newcommand{\blk}{$\upshape{blk}$}

\begin{document}
%
\title{Resilient distributed resource allocation algorithm under false data injection attacks}
%
%
%
\author{\IEEEauthorblockN{Xin~Cai, Xinyuan~Nan, Binpeng~Gao}

\thanks{
This work was supported by the National Natural Science Foundation of China (NSFC, Grant No. 52065064) and by the Natural Science Foundation of Xinjiang Province (Grant No. 2019D01C079). (Corresponding author: Xin~Cai.)

X. Cai, X. Nan and B. Gao are with the School of Electrical Engineering, Xinjiang University, Urumqi 830047, China (emails: xincai@xju.edu.cn; xynan@xju.edu.cn; gbp\_xd@sina.com). }}

\markboth{Transactions on }
{Cai \MakeLowercase{\textit{et al.}}:Resilient distributed resource allocation algorithm under false data injection attacks}

%



\maketitle

\begin{abstract}
 A resilient distributed algorithm is proposed to solve the distributed resource allocation problem of a first-order nonlinear multi-agent system who is subject to false data injection (FDI) attacks. An intelligent attacker injects false data into agents' actuators and sensors such that agents execute the algorithm according to the compromised control inputs and interactive information. The goal of the attacker is to make the multi-agent system to be unstable and to cause the deviance of agents' decisions from the optimal resource allocation. At first, we analyze the robustness of a distributed resource allocation algorithm under FDI attacks. Then, the unknown nonlinear term and the false data injected in agents are considered as extended states which can be estimated by extended state observers. The estimation was used in the feedback control to suppress the effect of the FDI attacks. A resilient distributed resource allocation algorithm based on the extended state observer is proposed to ensure that it can converge to the optimal allocation without requiring any information about the nature of the attacker. An example is given to illustrate the results.
\end{abstract}


%
\IEEEpeerreviewmaketitle

\section{Introduction}
As a special distributed constrained optimization problem, distributed resource allocation problem has attacked a lot of researchers' attention in various engineering applications, such as the economic dispatch in power systems \cite{Yi.2016,Pesaran.2017,Wang.2019}, the congestion control in data or traffic networks \cite{Low.1999,Farhadi.2019}, and the power optimization in sensor networks \cite{Lin.2006,Bishop.2010}. During the past decade, many continuous-time distributed algorithms have been designed for agents to obtain the optimal allocation only by local information.

Generally speaking, the computation and the communication have been  two important components in  the distributed optimization algorithms for solving the resource allocation problem. In recent years, the results of asymptotic and exponential convergence of distributed resource allocation algorithms have been obtained for agents communicating on undirected graphs\cite{Yi.2016,Li.2020}. And the related results were also extended to weight-balanced directed graphs (digraphs)\cite{Liang.2018,Zhu.2019,Chen.2021}. Furthermore, a distributed algorithm was designed for agents communicating on weight-unbalanced digraphs\cite{Zhu.2019}. In addition to the asymptotic and exponential convergence results, recent work dealt with finite-time and predefine-time convergence of distributed algorithms\cite{Guo.2020}. To cope with the more general cases, such as nonsmooth cost functions and heterogeneous local constraints, distributed optimization algorithms combined with differential inclusions and projection dynamics were proposed in \cite{Deng.2020c}. Additionally, a robust distributed resource allocation algorithm was designed to deal with the case of uncertain allocation parameters \cite{Zeng.2018}. Considering the problem with unknown cost functions, distributed algorithms based on the extremum seeking control were designed in \cite{Wang.2019,Ogwuru.2020}. The communication delays, which may exist in networks, were also considered in the implementation of distributed algorithms\cite{WangXF.2019}. To reduce the communication load, a distributed event-triggered algorithm was proposed in \cite{Deng.2020b}.  Recently, time-varying graphs incurred by persistent cyber attacks were considered in the resource allocation problem\cite{Wang.2020}. The abovementioned literatures focus on the problem with a fixed optimal allocation. Wang \textit{et al.} studied the problem with time-varying cost functions and designed distributed algorithms based on the prediction-correction method and nonsmooth consensus idea to converge to the time-varying optimal allocation\cite{Wang.2021,Wang.2020}. In addition, Deng \textit{et al.} studied the problem with the double-integrator and the multi-integrator agents and designed distributed algorithms based on the state feedback\cite{Deng.2018,Deng.2020,Deng.2020b,Deng.2020c}.  From the above observation, there is still much to be studied in the problem of distributed resource allocation, such as packet loss and cyber attacks in communication networks, and agents with nonlinear or more complex dynamics.

Note that most of the existing algorithms were designed under the setup of reliable communication networks. However, there may exist attackers who are latent in cyber systems to transmit false information in the network to destroy the control target of the multi-agent system. This type of cyber attacks is called false data injection (FDI) attacks. Besides, the denial of service (DoS) attacks, as another type of cyber attacks, also have attracted much attention in the filed of cyber-physical systems. Recently, the stability of a distributed resource allocation algorithm under DoS attacks was explored by the method of switched algorithm modeled as a hybrid system\cite{Shao.2020}. To the best of our knowledge, there is little work on the robustness of distributed resource allocation algorithms under FDI attacks. Here, we consider that an intelligent attacker lurks in a multi-agent system who aims to execute the designed algorithm to obtain the optimal allocation. The attacker injects false data into agents' actuators and sensors to manipulate the data in control inputs and the information obtained from communication networks. His purpose is to make the important data of physical systems wrong such that agents fail to obtain the optimal resource allocation. It is obvious that the attacker can easily succeed if no protective measurements are taken. To prevent the attacker to destroy the performance of the distributed algorithm, a resilient distributed resource allocation algorithm is designed based on the extended state observers. In summary, the main contributions of this paper are given below.

1) A distributed resource allocation problem under FDI attacks is formulated in this paper. The FDI attacks incur the failure of the distributed resource allocation algorithm by injecting false data into agents' actuators and sensors. Moreover, the attacks considered here are not detected and removed by identification.

2) The robustness of the distributed resource allocation algorithm for a multi-agent system with first-order nonlinear dynamics under FDI attacks is analyzed. The convergent condition is given to ensure that agents' decisions can be steered to a neighborhood of the optimal allocation only if agents suffer from weak FDI attacks.

3) A resilient distributed resource allocation algorithm is proposed for agents to obtain the optimal allocation under FDI attacks. We consider the nonlinear term in agents' dynamics and false data injected in the actuators and sensors as extended states. Extend state observers are used to restrain the influence of FDI attacks on the designed algorithm. Moreover, the designed resilient algorithm can deal with both FDI attacks and nonlinear (or uncertainties) existing in agents' dynamics.

The rest of this paper is organized as follows. In Section \uppercase\expandafter{\romannumeral2}, the problem formulation is given. In Section \uppercase\expandafter{\romannumeral3}, the robustness of a continuous-time distributed resource allocation algorithm is analyzed. A resilient distributed algorithm is designed in Section \uppercase\expandafter{\romannumeral4}. A simulation example is provided in Section \uppercase\expandafter{\romannumeral5}. Finally, some conclusions and future topics are stated in Section \uppercase\expandafter{\romannumeral6}.

Notations: $\mathbb{R}$ and $\mathbb{R}_{\geq0}$ denote the set of real and non-negative real numbers, respectively. $\mathbb{R}^n$ is the $n$-dimensional real vector space. $\mathbb{R}^{n\times m}$ denotes the set of $n \times m$ real matrices. Given a vector $x\in\mathbb{R}^n$, $\|x\|$ is the Euclidean norm. $A^T$ and $\|A\|$ are the transpose and the spectral norm of matrix $A\in \mathbb{R}^{n\times n}$, respectively. For matrices $A$ and $B$, $A\otimes B$ denotes their Kronecker product. Let $\col(x_1,\ldots,x_n)=[x_1^T,\ldots,x_n^T]^T$. Given matrices $A_1,\ldots,A_n$, $\blk\{A_1,\ldots,A_n\}$ denotes the block diagonal matrix with $A_i$ on the diagonal. $1_n$ and $0_n$ are $n$-dimensional column vectors with all elements being ones and zeros, respectively. $I_n$ denotes an $n\times n$ identity matrix.

\section{Problem formulation}
In this section, the basic problem of distributed resource allocation, a distributed resource allocation algorithm under reliable communication network, the model of FDI attacks and the algorithm under FDI attacks are given.

\subsection{Resource Allocation Problem}
A network of agents indexed in the set $\mathcal{I}=\{1,\ldots,N\}$ cooperates with each other to obtain an optimal allocation of the limited network resource, which can be formulated as following problem \eqref{op}.

\begin{equation} \label{op}
 \begin{array}{l}
\min f(X)=\sum_{i=1}^Nf_i(x_i)\\
s.t. \ \ \sum_{i=1}^Nx_i=\sum_{i=1}^Nd_i.
 \end{array}
\end{equation}
In \eqref{op}, $x_i\in\mathbb{R}^n$, $f_i(x_i):\mathbb{R}^n\rightarrow \mathbb{R}$, and $d_i\in\mathbb{R}^n$ are the decision variable, cost function and accessible resource data of agent $i$, respectively. $X=\col(x_1,\ldots,x_N)$ is the resource allocation vector of the whole network.
Distributed resource allocation problem is that agents in the network cooperatively find an optimal allocation according to their local information including costs, resource data and shared information with neighbors.

In the following, some basic assumptions are given and they are widely used in \cite{Yi.2016,Deng.2018,Zhu.2019,Li.2020,Deng.2020}.

\begin{assumption} \label{as1}
For each $i\in\mathcal{I}$, $f_i(x_i)$ is a continuously differential and strongly convex function with $l$-Lipschitz continuous gradient.
\end{assumption}

\begin{assumption} \label{as2}
There exists a finite optimal allocation $X^*$ for problem \eqref{op}.
\end{assumption}

\subsection{Agent's Dynamics}

We consider that each agent in the network has inherent dynamics, which are modeled by the following first-order nonlinear system
\begin{equation} \label{ad}
\dot{x}_i=g(x_i)+u_i, \ i\in\mathcal{I}.
\end{equation}
In \eqref{ad}, $x_i\in\mathbb{R}^{n}$ and $u_i\in\mathbb{R}^{n}$ are the decision and control input of agent $i$, respectively. $g_i(x_i):\mathbb{R}^n\rightarrow \mathbb{R}^n$ is a nonlinear function. The explicit expression of $g_i$ may be unknown.  For agent $i$, its control input is designed to steer its decision $x_i$ to reach the optimal allocation $x_i^*$ belonging to $X^*$.

To realize a distributed setting, each agent needs to share some information with its neighbors on a fixed undirected and connected graph, which is denoted by $\mathcal{G}$ $=$ $(\mathcal{I},\mathcal{E})$. $\mathcal{I}$ is the set of nodes corresponding to agents in the multi-agent system and $\mathcal{E}$ is the set of edges corresponding to communication links between neighboring agents. The Laplacian matrix of graph $\mathcal{G}$ is denoted by $L=[l_{ij}]\in\mathbb{R}^{N\times N}$, where $l_{ii}=\sum_{j=1}^N a_{ij}$ and $l_{ij}=-a_{ij}, i\neq j$ with the weight $a_{ij}$ on edge $(i,j)\in\mathcal{E}$. The eigenvalues of $L$ can be denoted by $0=\rho_1<\rho_2\leq \ldots \leq \rho_N$. For agent $i$, the set composed by its neighbors is denoted by $\mathcal{N}_i$. The property of $L$ is that $1_N^TL=L1_N=0_N$ \cite{Godsil.2001}.

\subsection{Distributed Resource Allocation Algorithm}
For agent $i$ with the known nonlinear term $g_i(x_i)$, the control input $u_i$ in \eqref{ad} is designed to steer its strategy $x_i$ to the optimal allocation. Thus,
\begin{equation} \label{alg1}
u_i=-g(x_i)-\nabla f_i(x_i)-\lambda_i, \ \ \forall i\in\mathcal{I},
\end{equation}
where $\nabla f_i(x_i)$ is the gradient of $f_i$ with respect to $x_i$, and $\lambda_i\in\mathbb{R}^n$ is the estimation of the Lagrange multiplier associated with equality constraint $\sum_{j=1}^Nx_j=\sum_{j=1}^Nd_j$. In \eqref{alg1}, $\lambda_i$ is regulated by the following dynamics,
\begin{equation} \label{lamda}
\begin{aligned}
\dot{\lambda}_i&=-\sum_{j=1}^N a_{ij}(\lambda_i-\lambda_j)-\sum_{j=1}^N a_{ij}(z_i-z_j)+x_i-d_i,\\
\dot{z}_i&=\sum_{i=1}^N a_{ij}(\lambda_i-\lambda_j),
\end{aligned}
\end{equation}
where $z_i\in\mathbb{R}^n$ is the auxiliary state that facilitate $\lambda_i$ to reach consensus.

\subsection{Attack Model}
\begin{figure}[!t]
  \centering
  \includegraphics[width=8cm]{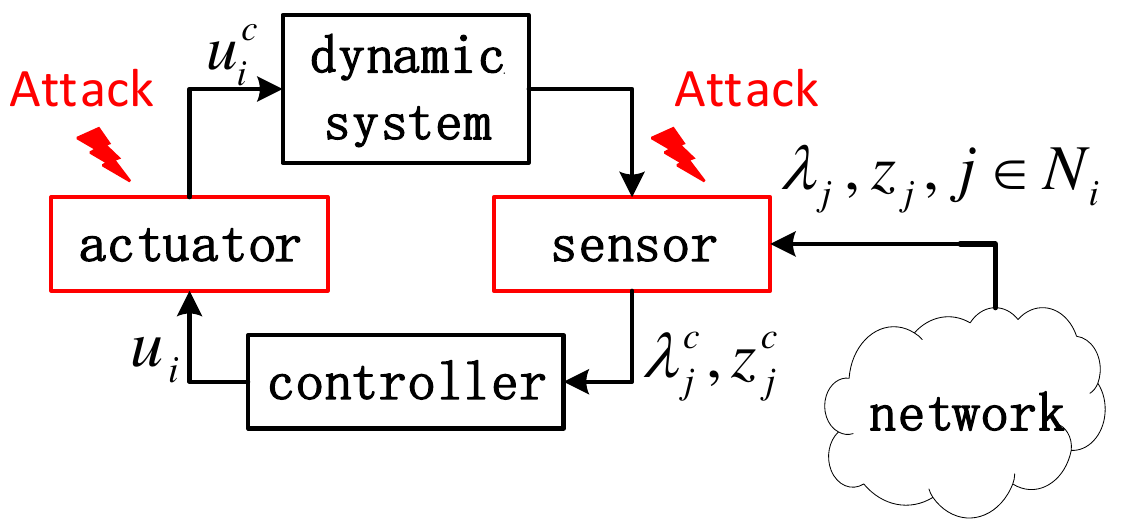}
  \caption{False data injection attacks on agent $i\in\mathcal{I}$ in a multi-agent system}
  \label{fig5}
\end{figure}

In this paper, an intelligent attacker injects false data into agents' actuators and sensors to prevent the multi-agent systems to find the optimal allocation of problem \eqref{op}. As shown in Fig.~\ref{fig5}, for agent $i$, the compromised input $u_i^c(t)$ received by the actuator can be modeled by
\begin{equation} \label{au}
u^c_i(t)=u_i(t)+u^a_i(t),
\end{equation}
where $u_i(t)$ is the uncompromised control input given in \eqref{alg1}, $u^a_i(t)$ is the unknown attack signal. The information sensed from neighboring agent $j~(j\in\mathcal{N}_i)$ can be modeled by
\begin{equation} \label{ss1}
\lambda_{j}^c(t)=\lambda_j(t)+\lambda_j^a(t),
\end{equation}
\begin{equation} \label{ss2}
z_j^c(t)=z_j(t)+z_j^a(t),
\end{equation}
where $\lambda_j(t)$ and $z_j(t)$ are the uncompromised information transmitted to agent $i$, $\lambda_j^{a}(t)$ and $z_j^a(t)$ are unknown attacks on the sensor of agent $i$, and $\lambda_j^{c}(t)$ and $z_j^c(t)$ are the compromised information sensed by agent $i$.

\subsection{Distributed Resource Allocation Algorithm under FDI Attacks}
For agent $i$ under FDI attacks modeled by \eqref{au}-\eqref{ss2}, the compromised distributed resource allocation algorithm is given by
\begin{equation}
\begin{aligned}
\dot{x}_i&=u_i^c,\\
\dot{\lambda}_i&=-\sum_{j=1}^N a_{ij}(\lambda_i^c-\lambda_j^c)-\sum_{j=1}^N a_{ij}(z_i^c-z_j^c)+x_i-d_i,\\
\dot{z}_i&=\sum_{i=1}^N a_{ij}(\lambda_i^c-\lambda_j^c).
\end{aligned}
\end{equation}
Sometime $t$ is omitted for simplicity. In the following development, for notional convenience, let $\kappa_{i1}=u_i^a$, $\kappa_{i2}=-\sum_{j=1}^N a_{ij}(\lambda_i^a-\lambda_j^a)-\sum_{j=1}^N a_{ij}(z_i^a-z_j^a)$, and $\kappa_{i3}=\sum_{j=1}^N a_{ij}(\lambda_i^a-\lambda_j^a)$. The distributed resource allocation algorithm under FDI attacks can be rewritten as follows.
\begin{equation} \label{alg2}
\begin{aligned}
\dot{x}_i&=-\nabla f_i(x_i)-\lambda_i+\kappa_{i1},\\
\dot{\lambda}_i&=-\sum_{j=1}^N a_{ij}(\lambda_i-\lambda_j)-\sum_{j=1}^N a_{ij}(z_i-z_j)+x_i-d_i+\kappa_{i2},\\
\dot{z}_i&=\sum_{i=1}^N a_{ij}(\lambda_i-\lambda_j)+\kappa_{i3}.
\end{aligned}
\end{equation}

Here, we only consider that the attacker has limited power. So, the boundedness of FDI attacks is assumed below.
\begin{assumption} \label{as3}
The FDI attacks $u_i^a(t)$, $\lambda_i^a(t)$, and $z_i^a(t)$ are bounded for all $i\in\mathcal{I}$. And $\dot{u}_i^a(t)$, $\dot{\lambda}_i^a(t)$, and $\dot{z}_i^a(t)$ are also bounded.
\end{assumption}

\begin{remark}
The boundedness of attacks indicates not only the stealthy property but also limited ability to falsify data. To the best of our knowledge, the FDI attacks satisfying Assumption \ref{as3} include uniformly bounded attacks\cite{Jin.2017,Yucelen.2016,Meng.2020}, attacks generated by exogenous systems\cite{Gusrialdi.2014,Gusrialdi.2018,Mustafa.2020,Huang.2020}, state-dependent attacks\cite{Jin.2017,Yucelen.2016,Dong.2020,Modares.2020}, and so on.
\end{remark}

\section{Robustness of Distributed Resource Allocation Algorithm}
In this section, the convergence of distributed resource allocation algorithm \eqref{ad}-\eqref{lamda} is given. Then, the robustness of distributed algorithm \eqref{ad}-\eqref{lamda} against FDI attacks is analyzed.
\subsection{Convergence}
Here, we consider that the communication network is secure and reliable. The convergence of distributed algorithm \eqref{ad}-\eqref{lamda} is analyzed by studying the stability of the closed-loop system which is given by
\begin{equation} \label{alg}
\begin{aligned}
\dot{x}_i&=-\nabla f_i(x_i)-\lambda_i,\\
\dot{\lambda}_i&=-\sum_{j=1}^N a_{ij}(\lambda_i-\lambda_j)-\sum_{j=1}^N a_{ij}(z_i-z_j)+x_i-d_i,\\
\dot{z}_i&=\sum_{i=1}^N a_{ij}(\lambda_i-\lambda_j).
\end{aligned}
\end{equation}
System \eqref{alg} is similar to the differentiated projected algorithm without considering the local convex constraints proposed in \cite{Yi.2016}. The relationship between the optimal allocation and the equilibrium of dynamic system \eqref{alg} is given in the following Lemma \ref{lemma1}. Let $X=\col(x_1,\ldots,x_N)$, $\Lambda=\col(\lambda_1,\ldots,\lambda_N)$, $Z=\col(z_1,\ldots,z_N)$, and $d=\col(d_1,\ldots,d_N)$. The overall system composed of $N$ agents with algorithm \eqref{alg} is described by
\begin{equation} \label{cf1}
\begin{aligned}
\dot{X}&=-F(X)-\Lambda,\\
\dot{\Lambda}&=-(L\otimes I_n)\Lambda-(L\otimes I_n)Z+X-d,\\
\dot{Z}&=(L\otimes I_n)\Lambda.
\end{aligned}
\end{equation}

\begin{lemma} \label{lemma1}
Under Assumptions \ref{as1}-\ref{as2}, $X^*$ is the optimal allocation of problem \eqref{op} if and only if there exist $1_N\otimes \lambda^*$ and $Z^*$ such that $(X^*, 1_N\otimes \lambda^*, Z^*)$ is the equilibrium point of system \eqref{cf1}.
\end{lemma}

Proof:  If $(X^*, \Lambda^*, Z^*)$ with $\Lambda^*=1_N\otimes \lambda^*$ is the equilibrium point of system \eqref{cf1}, we have that
\begin{subequations}
\begin{align}
0_{Nn}&=-F(X^*)-\Lambda^*,\label{cf11} \\
0_{Nn}&=-(L\otimes I_n)\Lambda^*-(L\otimes I_n)Z^*+X^*-d,\label{cf12} \\
0_{Nn}&=(L\otimes I_n)\Lambda^*. \label{cf13}
\end{align}
\end{subequations}
For an undirected and connected graph $\mathcal{G}$, it follows from \eqref{cf13} that $\lambda_i^*=\lambda_j^*=\lambda^*$. Left-multiplying \eqref{cf12} by $(1_N^T\otimes I_n)$, we have that
\begin{equation}\label{ec1}
\sum_{i=1}^N x_i^*-\sum_{i=1}^Nd_i=0_n.
\end{equation}
Furthermore, it follows from \eqref{cf11} that
\begin{equation} \label{ec2}
\nabla f_i(x_i^*)=\nabla f_j(x_j^*)=-\lambda^*, \ \forall i,j\in\mathcal{I}.
\end{equation}

Under Assumptions \ref{as1}-\ref{as2}, resource allocation problem \eqref{op} has a unique primal-dual solution $(x_i^*,\lambda^*)$ satisfying the following KKT conditions \cite{Xiao.2006}.
\begin{equation} \label{oc1}
\nabla f_i(x_i^*)+\lambda^*=0,
\end{equation}
and
\begin{equation} \label{oc2}
\sum_{j=1}^N x_j^*-\sum_{j=1}^N d_j=0.
\end{equation}
It follows from \eqref{ec1},\eqref{ec2},\eqref{oc1}, and \eqref{oc2} that the equilibrium point of system \eqref{cf1} is the optimal allocation of problem \eqref{op}, and vice versa. $\hfill\blacksquare$

The following lemma gives the result of the convergence of algorithm \eqref{alg}, which has been analyzed in \cite{Yi.2016}. We omit the proof of Lemma \ref{lemma2} for the limited space.

\begin{lemma}[{\cite[Theorem 4.3]{Yi.2016}}] \label{lemma2}
Suppose that Assumptions \ref{as1}-\ref{as2} hold. The distributed resource allocation algorithm \eqref{alg} can exponentially converge to the optimal allocation $X^*$ of problem \eqref{op}.
\end{lemma}

\subsection{Robustness}
Next, the robustness of distributed algorithm \eqref{alg} under FDI attacks is given in the following Theorem \ref{them2}.

\begin{theorem} \label{them2}
Suppose that Assumptions \ref{as1}-\ref{as3} hold. The distributed resource allocation algorithm \eqref{alg} under FDI attacks can converge to a neighborhood of the optimal allocation $X^*$ of problem \eqref{op}. The size of the neighborhood depends on the upper bound of the FDI attacks.
\end{theorem}

Proof: Distributed algorithm \eqref{alg} under FDI attacks can be described by \eqref{alg2}. The robustness of \eqref{alg} against the attacks is analyzed by studying the convergence of \eqref{alg2}. Denote $\tilde{X}=X-X^*$, $\tilde{\Lambda}=\Lambda-\Lambda^*$, $\tilde{Z}=Z-Z^*$, $\kappa_1=\col(\kappa_{11},\ldots,\kappa_{N1})$, $\kappa_2=\col(\kappa_{12},\ldots,\kappa_{N2})$, and $\kappa_3=\col(\kappa_{13},\ldots,\kappa_{N3})$.
The error system is
\begin{equation} \label{ce2}
\begin{aligned}
\dot{\tilde{X}}&=-(F(X)-F(X^*))-\tilde{\Lambda}+\kappa_1,\\
\dot{\tilde{\Lambda}}&=-(L\otimes I_N)\tilde{\Lambda}-(L\otimes I_N)\tilde{Z}+\tilde{X}+\kappa_2,\\
\dot{\tilde{Z}}&=(L\otimes I_N)\tilde{\Lambda}+\kappa_3.
\end{aligned}
\end{equation}
Let $\Theta=\col(\tilde{X},\tilde{\Lambda},\tilde{Z})$, $\nu(t)=\col(\kappa_1,\kappa_2,\kappa_3)$, and
\begin{align*}
\Gamma(\Theta)=\begin{bmatrix} -(F(X)-F(X^*))-\tilde{\Lambda}\\
-(L\otimes I_N)\tilde{\Lambda}-(L\otimes I_N)\tilde{Z}+\tilde{X}\\
(L\otimes I_N)\tilde{\Lambda}
\end{bmatrix}.
\end{align*}
Then, system \eqref{ce2} can be rewritten as
\begin{equation} \label{ps}
\dot{\Theta}=\Gamma(\Theta)+\nu(t),
\end{equation} which is a perturbed system with the nonvanishing perturbation $\nu(t)$. It follows from Lemma \ref{lemma2} that $\Theta^*=0_{3Nn}$ is an exponentially stable equilibrium point of the nominal system $\dot{\Theta}=\Gamma(\Theta)$. By the converse Lyapunov Theorem in \cite[Theorem 4.14]{Khalil.2002}, there is a continuous differential function $V:\mathbb{R}^{3Nn}\rightarrow \mathbb{R}$ that satisfies the inequalities
\begin{align*}
\begin{split}
c_1\|\Theta\|^2&\leq V \leq c_2\|\Theta\|^2,\\
\frac{\partial V}{\partial \Theta} \Gamma(\Theta)&\leq -c_3\|\Theta\|^2,\\
\|\frac{\partial V}{\partial \Theta}\|&\leq c_4\|\Theta\|
\end{split}
\end{align*}
for some positive constants $c_1$, $c_2$, $c_3$, and $c_4$. Then, the derivative of $V$ along the trajectories of perturbed system \eqref{ps} satisfies
\begin{align*}
\dot{V}&\leq -c_3\|\Theta\|^2+\|\frac{\partial V}{\partial \Theta}\|\|\nu(t)\|,\\
&\leq -c_3\|\Theta\|^2+c_4D\|\Theta\|\\
&\leq (1-\epsilon)\|\Theta\|^2, \ \forall \|\Theta\|\geq c_4D/\epsilon c_3,
\end{align*}
where $0<\epsilon<1$, the second inequality comes from the boundedness of the FDI attacks, that is, $\|\nu(t)\|\leq D$ for some $D>0$.
It further yields that $\lim_{t\rightarrow \infty}\|\Theta(t)\|<\sqrt{\frac{c_2}{c_1}}\frac{c_4D}{c_3\epsilon}$. It indicates that the distributed algorithm \eqref{alg1} in the case of the communication network in the presence of FDI attacks converges to a neighborhood of the optimal allocation $X^*$ of problem \eqref{op}. The size of the neighborhood depends on the upper bound of FDI attacks.

 $\hfill\blacksquare$

\begin{remark}
The convergence of distributed algorithm \eqref{alg} under the FDI attacks is guaranteed only if the upper bound of the FDI attacks is within an appropriate range. It is easily seen that the strength of FDI attacks may be too large to destroy the stability of system \eqref{alg2}. That is, the distributed algorithm \eqref{alg} can not converge to the optimal allocation if the FDI attacks are strong (see Figs. 2-3 in Section V).
\end{remark}

\section{Resilient Distributed Resource Allocation Algorithm}
In this section, the nonlinear term $g_i(x_i)$ is assumed to be unknown in dynamics \eqref{ad} for each agent $i\in\mathcal{I}$. Assume that $d{g}_i/dt$ exists and is bounded. An extended state observer is designed for agent $i$ to estimate the nonlinear term $g_i(x_i)$ and the influence of FDI attacks according to the measurements of $x_i$, $\lambda_i$, and $z_i$.

Let $\gamma_i=\col(x_i, \lambda_i, z_i)$ and $\kappa_i=\col(\kappa_{i1}+g_i, \kappa_{i2}, \kappa_{i3})$. The estimations of $\gamma_i$ and $\kappa_i$ are obtained by the following observer.
\begin{equation} \label{of}
\begin{aligned}
\dot{\hat{\gamma}}_i&=\hat{\kappa}_i+u_{io}+a_1h_1(\gamma_i-\hat{\gamma}_i),\\
\dot{\hat{\kappa}}_i&=a_2h_2(\gamma_i-\hat{\gamma}_i),
\end{aligned}
\end{equation}
where $\hat{\gamma}_i=\col(\hat{x}_i, \hat{\lambda}_i, \hat{z}_i)$, $\hat{\kappa}_i=\col(\hat{\kappa}_{i1}+\hat{g}_i, \hat{\kappa}_{i2}, \hat{\kappa}_{i3})$, $\hat{g}_i$ is the estimation of $g_i$, $a_1, a_2>0$, $h_1$ and $h_2$ are linear or nonlinear functions to be designed, and
\begin{align*}
u_{io}\!=\!\!\begin{bmatrix}
-\nabla f_i(x_i)-\lambda_i-\hat{\kappa}_{i1}-\hat{g}_i\\
\!-\!\!\sum_{j=1}^N\! a_{ij}(\lambda_i\!-\!\lambda_j)\!-\!\!\sum_{j=1}^N \!a_{ij}(z_i\!-\!z_j)\!+\!x_i\!-\!d_i\!-\!\hat{\kappa}_{i2} \\
\sum_{j=1}^N a_{ij}(\lambda_i-\lambda_j)-\hat{\kappa}_{i3}
\end{bmatrix}\!.
\end{align*}

The resilient distributed resource allocation algorithm is designed by
\begin{equation} \label{alg3}
\begin{aligned}
\dot{x}_i&=g_i(x_i)-\nabla f_i(x_i)-\lambda_i+\kappa_{i1}-\hat{\kappa}_{i1}-\hat{g}_i(x_i),\\
\dot{\lambda}_i&=-\!\sum_{j=1}^N a_{ij}(\lambda_i\!-\!\lambda_j)\!-\!\sum_{j=1}^N a_{ij}(z_i\!-\!z_j)\!+\!x_i\!-\!d_i\!+\!\kappa_{i2}\!-\!\hat{\kappa}_{i2},\\
\dot{z}_i&=\sum_{i=1}^N a_{ij}(\lambda_i-\lambda_j)+\kappa_{i3}-\hat{\kappa}_{i3}.
\end{aligned}
\end{equation}

\subsection{Linear ESO Based Resilient Distributed Algorithm}
The functions $h_1$ and $h_2$ are chosen by linear functions, that is, $h_1(\gamma_i-\hat{\gamma}_i)=w_0(\gamma_i-\hat{\gamma}_i)$ and $h_2(\gamma_i-\hat{\gamma}_i)=w_0^2(\gamma_i-\hat{\gamma}_i)$ with a positive constant $w_0$. The result on the convergence of linear ESO \eqref{of} is given in the Lemma \ref{lemma3}.

\begin{lemma} \label{lemma3}
Suppose that Assumption \ref{as3} holds. Then, linear ESO \eqref{of} with linear functions $h_1(\gamma_i-\hat{\gamma}_i)=w_0(\gamma_i-\hat{\gamma}_i)$ and $h_2(\gamma_i-\hat{\gamma}_i)=w_0^2(\gamma_i-\hat{\gamma}_i)$ has an arbitrary small observation error, as $w_0\rightarrow \infty$. That is,
\begin{align*}
\limsup_{t\rightarrow \infty} \|\gamma_i(t)-\hat{\gamma}_i(t)\| &\leq O((\frac{1}{w_0})^2),\\
\limsup_{t\rightarrow \infty} \|\kappa_i(t)-\hat{\kappa}_i(t)\| &\leq O(\frac{1}{w_0}),\ \ \forall i\in \mathcal{I}.
\end{align*}
\end{lemma}

Proof: To analyze the convergence of linear observer \eqref{of}, let $\eta_{i1}=\gamma_i-\hat{\gamma}_i$ and $\eta_{i2}=\frac{1}{w_0}(\kappa_i-\hat{\kappa}_i)$. Denote that $\eta_i=\col(\eta_{i1}, \eta_{i2})$, $\boldsymbol{A}=A\otimes I_{3n}$ with $A=\Big [\begin{smallmatrix} -a_1 & 1\\
-a_2 & 0 \end{smallmatrix} \Big ]$ and $\boldsymbol{B}=B\otimes I_{3n}$ with $B=\Big [\begin{smallmatrix}0\\ 1 \end{smallmatrix}\Big ]$. It follows from \eqref{of} that
\begin{equation} \label{eta}
\dot{\eta}_i=w_0\boldsymbol{A}\eta_i+ \frac{1}{w_0^2}\boldsymbol{B}\dot{\kappa}_{i}.
\end{equation}
Select $a_1=2$ and $a_2=1$ in \eqref{of} such that matrix $A$ is Hurwitz. Solving \eqref{eta}, we have
\begin{equation*}
\eta_i(t)=e^{w_0\boldsymbol{A}t}(\eta_i(0)+\frac{\delta}{w_0^2}\boldsymbol{A}^{-1}\boldsymbol{B})-\frac{\delta}{w_0^2}\boldsymbol{A}^{-1}\boldsymbol{B},
\end{equation*}
where $\|\dot{\kappa}_i\|\leq \delta$ by Assumption \ref{as3}.
Then, there exists a class $\mathcal{KL}$ function $\beta_1$ such that
\begin{equation} \label{etat}
\|\eta_i(t)\|\leq \beta_1(\|\eta_i(0)\|,t)+O((\frac{1}{w_0})^2).
\end{equation}
Recall the definition of $\eta_{i1}$ and $\eta_{i2}$. We have that
\begin{align*}
\limsup_{t\rightarrow \infty}\|\gamma_i(t)-\hat{\gamma}_i\| &\leq O((\frac{1}{w_0})^2),\\
\limsup_{t\rightarrow \infty}\|\kappa_i(t)-\hat{\kappa}_i\| &\leq O(\frac{1}{w_0}).
\end{align*}
Thus, we have the conclusion in Lemma \ref{lemma3}.
$\hfill\blacksquare$

\begin{theorem} \label{them3}
Suppose that Assumptions \ref{as1}-\ref{as3} hold. The resilient distributed resource allocation algorithm \eqref{alg3} with linear ESO \eqref{of} can converge to a small neighborhood of optimal allocation $X^*$ of problem \eqref{op}. That is,
\begin{align*}
\limsup_{t\rightarrow \infty}\|X-X^*\|=O(\frac{1}{w_0}).
\end{align*}
\end{theorem}

Proof: According to the definitions of $\kappa_i$ and $\hat{\kappa}_i$, we know that $\eta_{i2}=\col(\eta_{i21},\eta_{i22},\eta_{i23})$ with $\eta_{i21}=\frac{1}{w_0}(\kappa_{i1}-\hat{\kappa}_{i1}+g_i-\hat{g}_i)$, $\eta_{i22}=\frac{1}{w_0}(\kappa_{i2}-\hat{\kappa}_{i2})$, and $\eta_{i23}=\frac{1}{w_0}(\kappa_{i3}-\hat{\kappa}_{i3})$. Let $\eta_{21}=\col(\eta_{121},\ldots,\eta_{N21})$, $\eta_{22}=\col(\eta_{122},\ldots,\eta_{N22})$, and $\eta_{23}=\col(\eta_{123},\ldots,\eta_{N23})$. The compact form of \eqref{alg3} is given by
\begin{equation} \label{cf3}
\begin{aligned}
\dot{\tilde{X}}&=-(F(X)-F(X^*))-\tilde{\Lambda}+w_0\eta_{21},\\
\dot{\tilde{\Lambda}}&=-(L\otimes I_N)\tilde{\Lambda}-(L\otimes I_N)\tilde{Z}+\tilde{X}+w_0\eta_{22},\\
\dot{\tilde{Z}}&=(L\otimes I_N)\tilde{\Lambda}+w_0\eta_{23}.
\end{aligned}
\end{equation}
Denote $\nu'(t)=\col(\eta_{21},\eta_{22},\eta_{23})$. System \eqref{cf3} can be rewritten as
\begin{equation} \label{ps1}
\dot{\Theta}=\Gamma{\Theta}+w_0\nu'(t),
\end{equation}
which is a perturbed system with nonvanishing perturbation $\nu'(t)$.
Similar to the analysis of system \eqref{ps}, the derivative of $V$ along the trajectories of perturbed system \eqref{ps1} satisfies
\begin{align*}
\dot{V}&\leq -c_3\|\Theta\|^2+w_0\|\frac{\partial V}{\partial \Theta}\|\|\nu'(t)\|,\\
&\leq -c_3\|\Theta\|^2+w_0c_4\|\Theta\|\|\nu'\|.
\end{align*}

It follows from Lemma \ref{lemma3} and the input-to-state stability \cite{Khalil.2002} that there exist a class $\mathcal{KL}$ function $\sigma_1$ and a constant $\zeta$ such that
\begin{equation*}
\begin{aligned}
\Theta_1(t)\leq \sigma_1(\Theta_1(0),t)+\zeta\sup_{0\leq \tau \leq t}\|\nu'(\tau)\|.
\end{aligned}
\end{equation*}
Denote $\Theta_2=\col(X-X^*, \Lambda-\Lambda^*, Z-Z^*, \eta)$. According to \eqref{etat} and the proof of Lemma 4.7 in \cite{Khalil.2002}, it yields that
\begin{equation*}
\Theta_2(t)\leq \sigma_2(\Theta(0),t)+O(\frac{1}{w_0}).
\end{equation*}
Thus, we have that $\limsup_{t\rightarrow \infty} \|X-X^*\|=O(\frac{1}{w_0})$. It indicates that the designed resilient algorithm can converge to a small neighborhood of optimal allocation $X^*$ of problem \eqref{op} even if agents are subject to the FDI attacks.
$\hfill\blacksquare$

\subsection{Nonlinear ESO Based Resilient Distributed Algorithm}
Since linear ESOs may cause overshot in the transient process of system \eqref{alg3}\cite{Khalil.2017}, nonlinear ESOs can be used to overcome this weakness by selecting appropriate nonlinear function $h_1$ and $h_2$ in \eqref{of}. Here, $h_1$ and $h_2$ are need to satisfy Assumption \ref{as4}, and can be selected by the sign function \cite{Guo.2011,Ren.2019}, or the piecewise function designed in \cite{Han.2009}.

\begin{assumption} \label{as4}
There exist constants $l_i (i=1,2,3,4,5)$ and positive definite, continuous differentiable functions $V_1$, $W_1: \mathbb{R}^{6n}\rightarrow \mathbb{R}$ such that

1) $l_1\|\eta_i\|^2\leq V_1(\eta_i) \leq l_2\|\eta_i\|^2, \ \ l_3\|\eta_i\|^2\leq W_1(\eta_i) \leq l_4\|\eta_i\|^2,$

2) $\frac{\partial V_1}{\partial \eta_{i1}}(\eta_{i2}-h_1(\eta_{i1}))-\frac{\partial V_2}{\partial \eta_{i2}} h_{2}(\eta_{i1})\leq -W_1(\eta_i)$,

3) $|\frac{\partial V_1}{\partial \eta_{i2}}| \leq l_5 \|\eta_i\|$,

\noindent where $\eta_i=\col(\eta_{i1}, \eta_{i2})$ defined in \eqref{eta}.
\end{assumption}

\begin{lemma} \label{lemma4}
Suppose that Assumptions \ref{as3}-\ref{as4} hold. Let $a_1=w_0$ and $a_2=w_0^2$.  Then, the nonlinear ESO \eqref{of} has an arbitrarily small observation error, as $w_0\rightarrow \infty$. That is,
\begin{align*}
\limsup_{t\rightarrow \infty} \|\gamma_i(t)-\hat{\gamma}_i(t)\| &\leq O((\frac{1}{w_0})^2),\\
\limsup_{t\rightarrow \infty} \|\kappa_i(t)-\hat{\kappa}_i(t)\| &\leq O(\frac{1}{w_0}),\ \ \forall i\in \mathcal{I}.
\end{align*}
\end{lemma}

Proof: According the definition of $\eta_i$ in \eqref{eta}, it follows from \eqref{of} that
\begin{equation*}
\begin{aligned}
\dot{\eta}_{i1}&=w_0(\eta_{i2}-h_1(\eta_{i1})),\\
\dot{\eta}_{i2}&=-w_0h_2(\eta_{i1})+\frac{1}{w_0}\dot{\kappa}_{i},
\end{aligned}
\end{equation*}
where $a_1=w_0$ and $a_2=w_0^2$.
Under Assumption \ref{as4}, we have that
\begin{align*}
\begin{split}
\dot{V}_1(\eta_i)&\!=\!\frac{\partial V_1}{\partial \eta_{i1}}(w_0(\eta_{i2}\!-\!h_1(\eta_{i1})))\\
&\ \ \ +\frac{\partial V_1}{\partial \eta_{i2}}(-w_0h_2(\eta_{i1})+\frac{1}{w_0}\dot{\kappa}_{i})\\
&\leq -w_0W_1(\eta_i)+\frac{\partial V_1}{\partial \eta_{i2}}\frac{1}{w_0}\dot{\kappa}_{i}.
\end{split}
\end{align*}
According to Assumption \ref{as4} and the definition of $\kappa_{i}$, we know that $\dot{\kappa}_{i}$ is bounded, that is, $\|\dot{\kappa}_{i}\|\leq \delta$. It is further deduced that
\begin{align*}
\dot{V}_1&\leq -w_0W_1(\eta_i)+\frac{l_5\delta}{w_0}\|\eta_i\|\\
& \leq -\frac{w_0l_{3}}{l_{2}}V_1(\eta_i)+\frac{l_5\delta\sqrt{l_{1}}}{l_{1}w_0}\sqrt{V_1(\eta_i)}.
\end{align*}
Thus, it yields that
\begin{align*}
\|\eta_i\|&\leq \frac{\sqrt{l_{1}}}{l_{1}}\sqrt{V_1}\\
&\leq \!\frac{\sqrt{l_{1}V_1(0)}}{l_{1}}e^{-\frac{w_0l_{3}}{2l_{2}}t}\!+\!\frac{l_5\delta}{2l_{1}w_0}\!\int_0^te^{-\frac{w_0l_{3}}{2l_{2}}(t-\tau)}d\tau\\
&\leq \frac{\sqrt{l_{1}V_1(0)}}{l_{1}}e^{-\frac{w_0l_{3}}{2l_{2}}t}+\frac{l_5\delta l_{2}}{l_{1}l_{3}w_0^2}(1-e^{-\frac{w_0l_{3}}{2l_{2}}t}).
\end{align*}
Then, $\limsup_{t\rightarrow \infty}\|\eta_i(t)\|\leq O(\frac{1}{w_0^2})$. It follows from the definition of $\eta_i$ that
\begin{align*}
\limsup_{t\rightarrow \infty}\|\gamma_i(t)-\hat{\gamma}_i\| &\leq O((\frac{1}{w_0})^2),\\
\limsup_{t\rightarrow \infty}\|\kappa_i(t)-\hat{\kappa}_i\| &\leq O(\frac{1}{w_0}).
\end{align*}
Thus, the conclusion in Lemma \ref{lemma4} is obtained. $\hfill\blacksquare$

\begin{remark}
It follows from Lemmas \ref{lemma3}-\ref{lemma4} that linear ESOs and nonlinear ESOs have the similar property of convergence. Whereas, compared with linear ESOs, nonlinear ESOs can make the transient process of the closed-loop system be with a smaller overshot by the feedback control (see Figs. 4-5 in Section V).
\end{remark}

\begin{theorem}
Suppose that Assumptions \ref{as1}-\ref{as4} hold. The designed resilient distributed resource allocation algorithm \eqref{alg3} with the nonlinear ESO \eqref{of} can converge to a small neighborhood of the optimal allocation $X^*$ of problem \eqref{op}. That is,
\begin{align*}
\limsup_{t\rightarrow \infty}\|X-X^*\|=O(\frac{1}{w_0}).
\end{align*}
\end{theorem}

The proof is similar to the proof of Theorem \ref{them3} and is omitted for the limited space.


\section{An example}

\begin{figure}[!t]
  \centering
  \includegraphics[width=8cm]{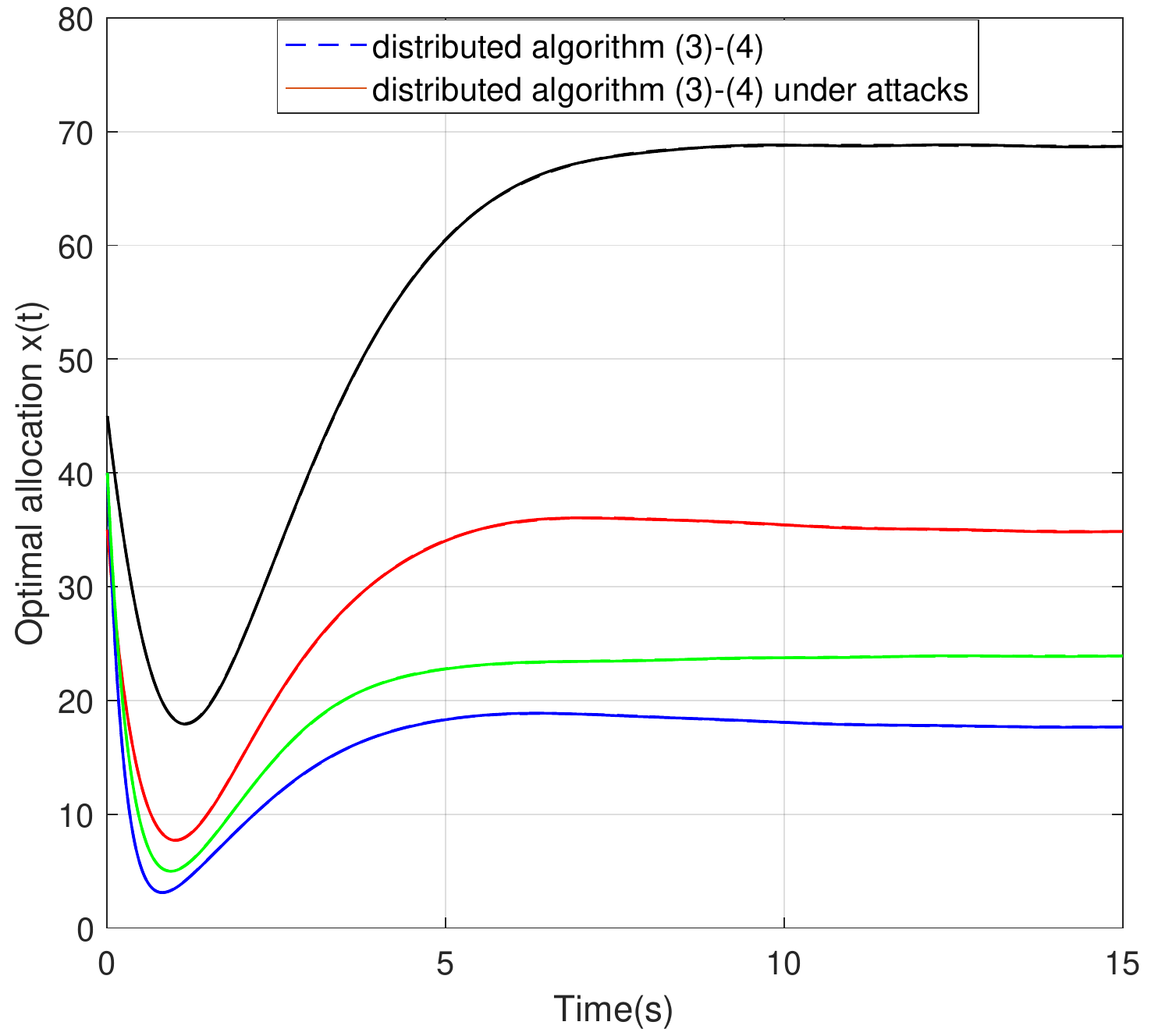}
  \caption{The optimal resource allocation under attacks with $u_i^a(t)=0.1\cos(2t)$, $\lambda^{ja}(t)=0.2\cos(2t)$, and $z^{ja}=0.1\cos(2t)$ $i\in\mathcal{I}, j\in\mathcal{N}_i$}
  \label{fig1}
\end{figure}

\begin{figure}[!t]
  \centering
  \includegraphics[width=8cm]{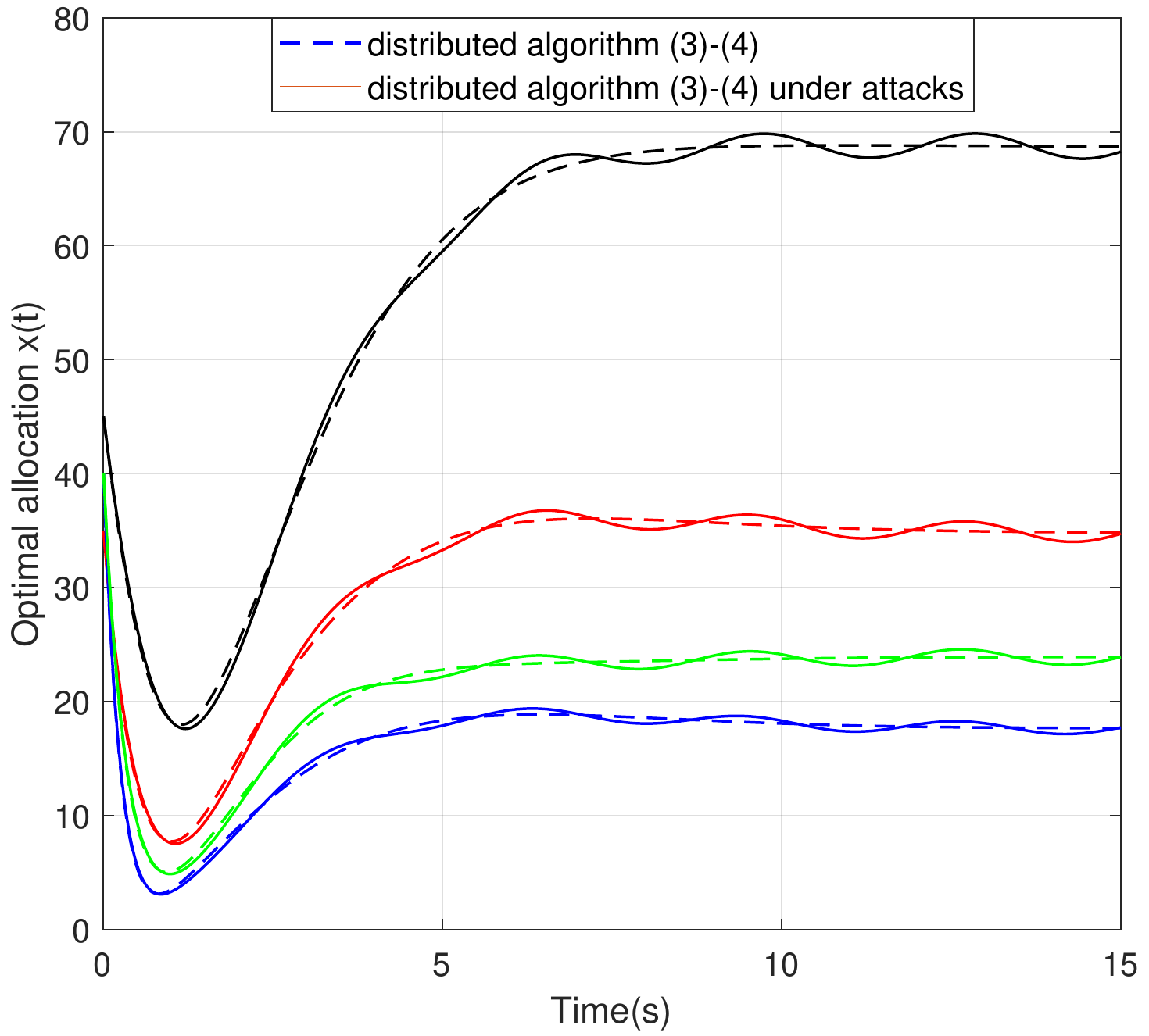}
  \caption{The optimal resource allocation under attacks with $u_i^a(t)=2\cos(2t)$, $\lambda^{ja}(t)=1.5\cos(2t)$, and $z^{ja}=\cos(2t)$ $i\in\mathcal{I}, j\in\mathcal{N}_i$}
  \label{fig2}
\end{figure}

\begin{figure}[!t]
  \centering
  \includegraphics[width=8cm]{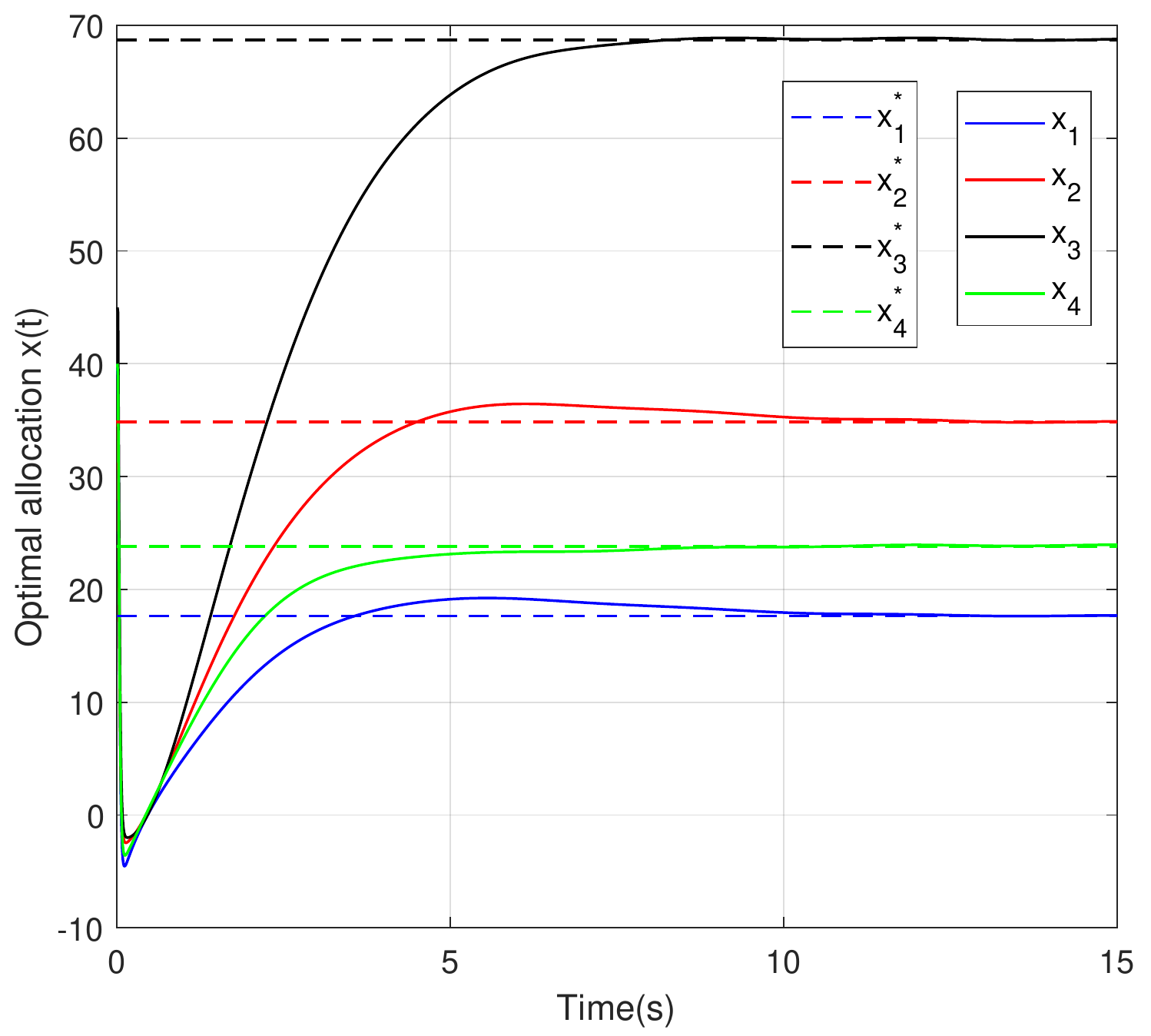}
  \caption{The optimal resource allocation obtained by the resilient distributed algorithm \eqref{alg3} based on linear ESO}
  \label{fig3}
\end{figure}

\begin{figure}[!t]
  \centering
  \includegraphics[width=8cm]{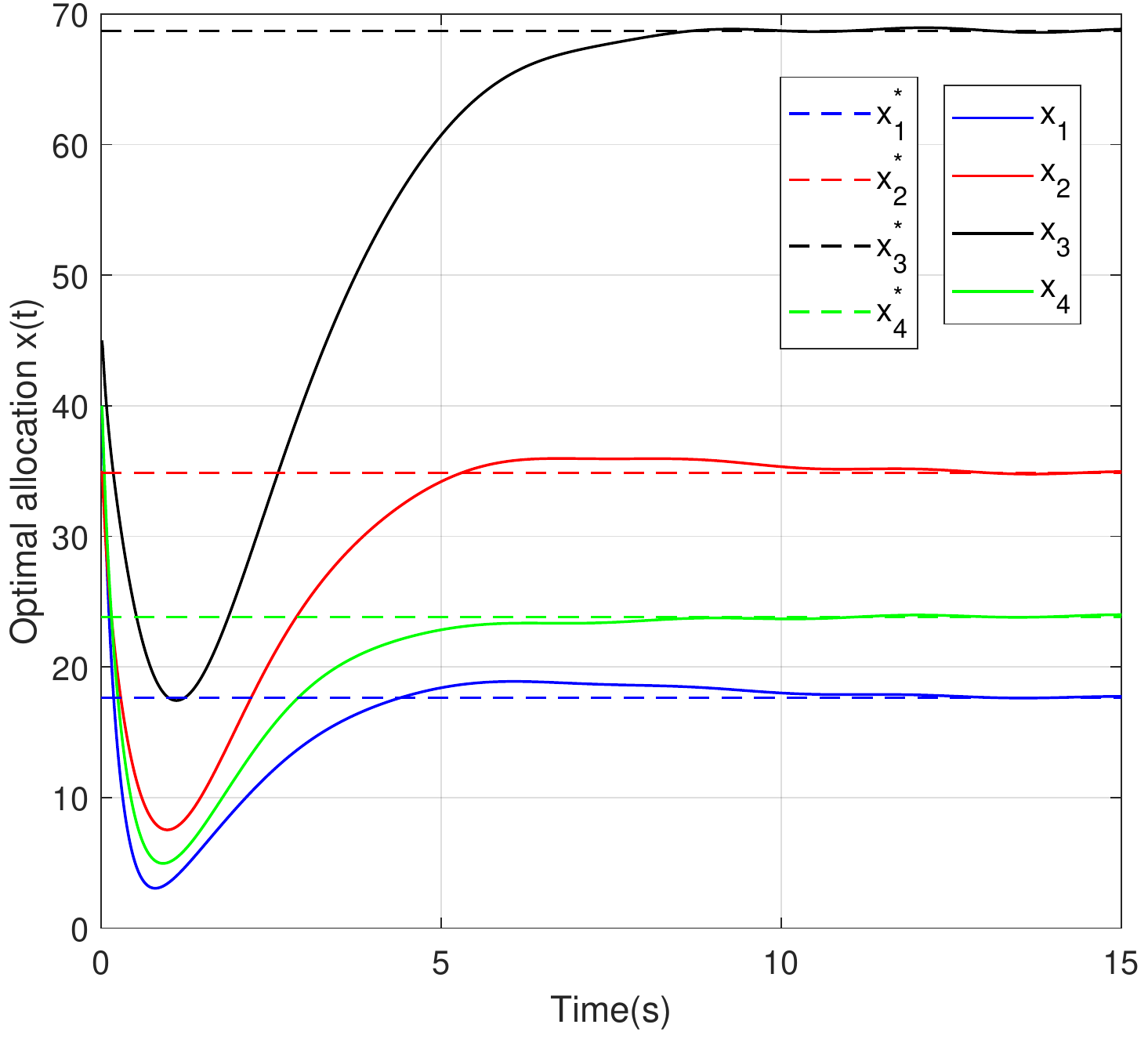}
  \caption{The optimal resource allocation obtained by the resilient distributed algorithm \eqref{alg3} based on nonlinear ESO}
  \label{fig4}
\end{figure}

In this section, we give an example of the economic dispatch problem of a power system with four generators to illustrate the obtained results. Each generator is considered as an agent and communicates with its neighbors on a line graph. For agent $i$, the nonlinear term in its dynamics \eqref{ad} is $g(x_i)=\sin(x_i)$, and the cost function is described by
\begin{equation*}
J_i(x_i,x_{-i})=a_i+b_ix_i+c_ix_i^2, \forall i\in\mathcal{I}.
\end{equation*}
Parameters in cost functions are given by $a=[a_1,\ldots,a_4]^T=[0.5,1.5,3.0,1.0]^T$, $b=[b_1,\ldots,b_4]^T=[3, 4, 5, 2]^T$ and $c=[c_1,\ldots,c_4]^T=[2, 1, 0.5, 1.5]^T$. The total demand of this power system is $d=145$, which is allocated by $[d_1,\ldots,d_4]^T=[30, 40, 40, 35]^T$. It is calculated that the optimal allocation is $X^*=[17.65,34.85,68.70,23.8]^T$. The initial state is $X(0)=[x_1(0),\ldots,x_4(0)]^T=[40, 35, 45, 40]^T$. In Fig.~\ref{fig1}, the FDI attacks are selected by $u_i^a(t)=0.1\cos(2t)$, $\lambda^{ja}(t)=0.2\cos(2t)$, and $z^{ja}=0.1\cos(2t)$. We can see that all agents' decisions reach the optimal resource allocation by the distributed algorithm \eqref{ad}-\eqref{lamda}. When this multi-agent system is subject to the FDI attacks with $u_i^a(t)=\cos(2t)$ $\lambda^{ja}(t)=1.5\cos(t)$, and $z^{ja}=\cos(2t)$ by increasing amplitudes of the attacks. The optimal allocation obtained by the algorithm \eqref{ad}-\eqref{lamda} under the attacks is shown in Fig.~\ref{fig2}. It is seen that all agents' decisions oscillate around the optimal allocation. From the comparison of Figs.~\ref{fig1}-\ref{fig2}, it indicates that the distributed algorithm \eqref{ad}-\eqref{lamda} is robust to the attacks with the smaller amplitudes.

In the resilient distributed resource allocation algorithm \eqref{alg3}, we select  $w=50$ in linear ESO \eqref{of}. The optimal allocation obtained by the linear ESO based resilient algorithm \eqref{alg3} is shown in Fig.~\ref{fig3}. It is seen that the strategies converge to optimal allocation $X^*$. However, the overshot phenomenon appears in the regulation process of strategies. Therefore, we apply the nonlinear ESO to replace the linear ESO by selecting $h_2=fal(\gamma_i-\hat{\gamma}_i,0.125,0.5)$ in \eqref{of}. The nonlinear function $fal(e,\alpha,\delta)$ defined by
\begin{equation*}
fal(e,\alpha,\delta)=
\begin{cases}
\frac{e}{\delta^{\alpha-1}} & |e|\leq \delta \\
|e|^\alpha sign(e) & |e|> \delta
\end{cases}.
\end{equation*}
The optimal allocation obtained by the resilient algorithm \eqref{alg3} based on the nonlinear ESO is shown in Fig.~\ref{fig4}. Compared with Fig.~\ref{fig3}, the transient process with the smaller overshot is shown in Fig.~\ref{fig4}.

\section{Conclusions}
In this paper, the robustness of a continuous-time distributed resource allocation algorithm under FDI attacks has been analyzed for a multi-agent system with first-order nonlinear dynamics. Then, to suppress the effect of the attacker's behavior, a resilient distributed algorithm was proposed based on the extended state observer, which was used to estimate the FDI attacks on agents. The sufficient condition for the convergence of the designed resilient algorithm was given by the Lyapunov stability theory and the stability of perturbed systems. In the future, it may be an interesting problem to design a resilient distributed resource allocation algorithm for a multi-agent system under DoS attacks.


%



\ifCLASSOPTIONcaptionsoff
  \newpage
\fi



%

\bibliographystyle{IEEEtran}
\bibliography{reference}

\end{document}